\documentclass[a4paper,12pt]{amsart}
\usepackage{mathrsfs}
\usepackage{amsfonts}
\usepackage{amssymb}
\usepackage{amsfonts, amscd, amsmath, mathrsfs, amssymb, amsthm, amsxtra, bbding, epsfig, graphicx, latexsym, url, mathbbol, bbold}
\usepackage[papersize={8in,11in},textwidth=15.7cm,textheight=21.3cm,centering]{geometry}
\usepackage{enumerate}

\usepackage[colorlinks=true,citecolor=blue,linkcolor=blue]{hyperref}
\hypersetup{
pdfstartpage=1,
pdfstartview=FitH}

\DeclareFontFamily{U}{mathx}{\hyphenchar\font45}
\DeclareFontShape{U}{mathx}{m}{n}{
      <5> <6> <7> <8> <9> <10>
      <10.95> <12> <14.4> <17.28> <20.74> <24.88>
      mathx10
      }{}
\DeclareSymbolFont{mathx}{U}{mathx}{m}{n}
\DeclareMathAccent{\widecheck}{\mathalpha}{mathx}{"71}
\DeclareFontFamily{U}{mathb}{\hyphenchar\font45}
\DeclareFontShape{U}{mathb}{m}{n}{
      <5> <6> <7> <8> <9> <10> gen * mathb
      <10.95> mathb10 <12> <14.4> <17.28> <20.74> <24.88> mathb12
}{}
\DeclareSymbolFont{mathb}{U}{mathb}{m}{n}
\DeclareMathSymbol{\llcurly}{3}{mathb}{"CE}
\DeclareMathSymbol{\ggcurly}{3}{mathb}{"CF}

\usepackage{caption} 
\numberwithin{equation}{section}

\allowdisplaybreaks

\newtheorem{theorem}{Theorem}[section]
\newtheorem{lemma}{Lemma}[section]

\makeatletter
\newcounter{roem}
\renewcommand{\theroem}{\Roman{roem}}

\newcommand{\c@org@eq}{}
\let\c@org@eq\c@equation
\newcommand{\org@theeq}{}
\let\org@theeq\theequation

\newcommand{\setroem}{
\let\c@equation\c@roem
 \let\theequation\theroem}

\newcommand{\setarab}{
\let\c@equation\c@org@eq
\let\theequation\org@theeq}
\makeatother

\newtheorem*{claim*}{Claim}

\theoremstyle{remark}
\newtheorem{remark}{\bf Remark}

\newcommand{\ue}{\mathrm{e}}
\newcommand{\ud}{\mathrm{d}}

\newcommand{\ft}{\mathrm{FT}}

\newcommand{\tr}{\mathrm{tr}}

\newcommand{\rank}{\mathrm{rank}}
\newcommand{\Swan}{\mathrm{Swan}}
\newcommand{\kl}{\mathrm{Kl}}

\DeclareMathOperator{\Mod}{mod}

\renewcommand{\bmod}[1]{\,(\Mod{ #1})}

\newcommand{\bk}{\mathbf{k}}
\newcommand{\bn}{\mathbf{n}}

\newcommand{\bA}{\mathbf{A}}
\newcommand{\bC}{\mathbf{C}}

\newcommand{\bF}{\mathbf{F}}

\newcommand{\bP}{\mathbf{P}}
\newcommand{\bQ}{\mathbf{Q}}
\newcommand{\bR}{\mathbf{R}}
\newcommand{\bZ}{\mathbf{Z}}

\newcommand{\cF}{\mathcal{F}}
\newcommand{\cG}{\mathcal{G}}

\newcommand{\cK}{\mathcal{K}}
\newcommand{\cL}{\mathcal{L}}

\newcommand{\fc}{\mathfrak{c}}

\newcommand{\fp}{\mathfrak{p}}
\newcommand{\fq}{\mathfrak{q}}

\newcommand{\fS}{\mathfrak{S}}

\newcommand{\sB}{\mathscr{B}}

\newcommand{\sD}{\mathscr{D}}

\newcommand{\sQ}{\mathscr{Q}}

\newcommand{\sS}{\mathscr{S}}

\def\leq{\leqslant}

\def\geq{\geqslant}
\def\llc{\llcurly}

\usepackage{graphicx}
\usepackage{tikz}

\begin{document}
\vglue -2mm

\title{A shifted convolution sum for $GL(3)\times GL(2)$}

\author{Ping Xi}

\address{Department of Mathematics, Xi'an Jiaotong University, Xi'an 710049, P. R. China}
\email{ping.xi@xjtu.edu.cn}




\begin{abstract}
In this paper, we estimate the shifted convolution sum
\[\sum_{n\geqslant1}\lambda_1(1,n)\lambda_2(n+h)V\Big(\frac{n}{X}\Big),\]
where $V$ is a smooth function with support in $[1,2]$, $1\leqslant|h|\leqslant X$, $\lambda_1(1,n)$ and $\lambda_2(n)$ are the $n$-th Fourier coefficients of $SL(3,\bZ)$ and $SL(2,\bZ)$ Hecke--Maass cusp forms, respectively. We prove an upper bound $O(X^{\frac{21}{22}+\varepsilon})$, updating a recent result of Munshi.
\end{abstract}

\maketitle


\section{Introduction}
\label{sec:Introduction}
Given two arithmetic sequences $\alpha$ and $\beta$, one usually encounters in analytic number theory the shifted convolution sum
\begin{align*}
\sum_{n\leqslant X}\alpha(n)\beta(n+h)
\end{align*}
with $h\neq0$. To state the status of this problem in principle, we would like to recall the following special cases that have been studied in literature:
\begin{itemize}
\item $(\alpha,\beta)=(\Lambda,\Lambda)$: this is related to the twin prime conjecture ($h=2$)
\item $(\alpha,\beta)=(\mu,\mu)$: this is related to the Chowla conjecture
\item $(\alpha,\beta)=(\mu^2,\mu^2)$: this counts consecutive squarefree numbers initiated by Heath-Brown $(h=1)$
\item $(\alpha,\beta)=(\Lambda,\tau)$: this is the focus of classical Titchmarsh divisor problem
\item $(\alpha,\beta)=(\tau,\tau)$: this is known as additive divisor problem, related to the fourth moments of Riemann zeta functions; one can also consider general cases with $\tau$ replaced by $\tau_k$
\item $(\alpha,\beta)=(\lambda_f,\lambda_g)$: this is the so-called shifted convolution problem for $GL(2)\times GL(2)$ as a cuspidal analogue of the additive divisor problem, which is quite important in the study of $GL(2)$ $L$-functions
\end{itemize}
Here $\Lambda,\mu,\tau$ denote the von Mangoldt, M\"obius and divisor functions, respectively; and $\lambda_f,\lambda_g$ denote the Fourier coefficients of $GL(2)$ cusp forms $f,g$. More details and relevant developments can be referred to \cite{Bl04,Bl05,BHM07,BFI86,DFI94,Ha03,HM06,HB84,Ho09,KMV02,LS03,MRT15,Mic04,Sa01,Se65} for instance.

In this paper, we consider the following smoothed version of shifted convolution sum
\begin{align*}
\sD_h(X)=\sum_{n\geqslant1}\lambda_1(1,n)\lambda_2(n+h)V\Big(\frac{n}{X}\Big),
\end{align*}
where $V$ is a fixed smooth function with compact support in $[1,2]$ and
\begin{itemize}
\item $\lambda_1(1,n)$ is the $n$-th Fourier coefficient of an $SL(3,\bZ)$ Hecke--Maass cusp form $\pi_1$,
\item $\lambda_2(n)$ is the $n$-th Fourier coefficient of an $SL(2,\bZ)$ Hecke--Maass or Hecke holomorphic cusp form $\pi_2$.
\end{itemize}
By Cauchy's inequality and Rankin--Selberg
theory (see Lemma \ref{lm:second-GL3} below), one has
\begin{align*}
\sD_h(X)\ll X^{1+\varepsilon},
\end{align*}
where the implied constant is allowed to depend on $\pi_1,\pi_2$ and of course on $\varepsilon.$ The interest of studying $\sD_h(X)$ is to capture the non-correlations between the $GL(3)$ and $GL(2)$ objects by obtaining a power saving in $X$.
Munshi \cite{Mu13a} is the first who beats the above (trivial) bound by saving an exponent $\frac{1}{26}$ in $X$ (after a suitable correction on his estimates for exponential sums in Lemma 11). The main ingredient in \cite{Mu13a} is a variant of the circle method due to Jutila \cite{Ju92,Ju97}, where the moduli can be chosen at one's demand. The innovation of Munshi is to specialize the moduli to be the product of two primes of certain sizes, so that he can utilize the bilinear structures in the estimates for resultant exponential sums. 
Following Munshi's approach, we may prove the following stronger power-saving by exploring the bilinear structure more efficiently in the circle method.

\begin{theorem}\label{thm:individual}
For $1\leqslant|h|\leqslant X,$ we have
\begin{align*}
\sD_h(X)\ll X^{\frac{21}{22}+\varepsilon},
\end{align*}
where the implied constant depends on $\varepsilon,\pi_1$ and $\pi_2.$
\end{theorem}

As noted in \cite{Mu13a}, the dependence of the implied constant on the conductors of $\pi_1,\pi_2$ can be made explicitly.

Much earlier than Munshi, Pitt \cite{Pi95} studied the similar shifted convolution sum with respect to the ternary divisor function $\tau_3(n)$ instead of $\lambda_1(1,n)$. One can also refer to \cite{Mu13b}, \cite{Ta15} and \cite{Su17} for recent progresses.

\subsection*{Notation and convention}
As usual, we write $\ue(t)=\ue^{2\pi it}$. The variable $p$ is reserved for prime numbers.
For a function $f$ defined over $\bZ/q\bZ,$ the Fourier transform is defined as
\begin{align*}
\widehat{f}(y) := \frac{1}{\sqrt{q}}\sum_{a\bmod{q}}f(a)\mathrm{e}\Big(-\frac{ya}{q}\Big).
\end{align*}
For a function $g\in L^1(\bR)$, its Fourier transform is defined as
\begin{align*}
\widehat{g}(y) := \int_\bR g(x) \mathrm{e}(-yx)\ud x.
\end{align*}

We use $\varepsilon$ to denote a very small positive number, which might be different at each occurrence; we also write $X^\varepsilon \log X\ll X^\varepsilon.$ The notation $n\sim N$ means $N<n\leqslant2N.$  
If $X$ and $Y$ are two quantities depending on $x$, we say that $X=O(Y)$ or $X\ll Y$ if one has $|X|\leqslant CY$ for some fixed $C$, and $X=o(Y)$
if one has $X/Y$ tends to zero as $x\rightarrow+\infty$. We use $X\llcurly Y$ to
denote the estimate $X\ll x^{o(1)}Y$.

\subsection*{Acknowledgements} 
The work is supported in part by NSFC (No.11601413).

\smallskip

\section{Background on automorphic forms}
\label{sec:Automorphicforms}

In this section, we recall some basic concepts and results on $GL(2)$ and $GL(3)$ automorphic forms.

We shall first briefly recall some basic facts about $SL(3,\bZ)$ automorphic forms, and the details can be found in Goldfeld's book \cite{Go06}. Suppose $\pi_1$ is a Maass form of type $(\nu_1,\nu_2)$ for $SL(3,\bZ)$ which is an eigenfunction of all the Hecke operators with Fourier coefficients $\lambda_1(m_1, m_2)$, normalized so that $\lambda_1(1, 1)=1$.  We introduce the Langlands parameters $(\alpha_1, \alpha_2, \alpha_3)$, defined by
\[
\alpha_1=-\nu_1-2\nu_2+1,\ \ \ 
\alpha_2=-\nu_1+\nu_2\ \ \ \text{and}\ \ \ 
\alpha_3=2\nu_1+\nu_2-1.
\]
The Ramanujan--Selberg conjecture predicts that $|\text{Re}(\alpha_i)|=0$, and from the work of Jacquet and Shalika we at least know that $|\text{Re}(\alpha_i)|<\frac{1}{2}$. 

We also recall the following (inverse) Hecke multiplicativity
\begin{align}\label{eq:Heckemultiplicativity}
\lambda_1(m_1,m_2)=\sum_{d|(m_1,m_2)}\mu(d)\lambda_1\Big(\frac{m_1}{d},1\Big)\lambda_1\Big(1,\frac{m_2}{d}\Big).
\end{align}

Let $g$ be a compactly supported function on $\bR^+$, and let 
$$\widetilde{g}(s)=\int_0^\infty g(x)x^{s-1}\ud x$$ be the Mellin transform. For $\sigma>-1+\max\{-\text{Re}(\alpha_1),-\text{Re}(\alpha_2),-\text{Re}(\alpha_3)\}$ and $\ell=0,1$ define
\[G_{\ell}(y)=\frac{1}{2\pi i}\int_{(\sigma)}\frac{\widetilde g(-s)}{(\pi^3 y)^{s}}\prod_{j=1}^3\frac{\Gamma(\frac{1+s+\alpha_j+\ell}{2})}{\Gamma(\frac{-s-\alpha_j+\ell}{2})}\ud s\]
and set
\begin{align}
\label{eq:Gpm}
G_\pm(y)=\frac{1}{2\pi^{3/2}}\left(G_{0}(y)\mp iG_{1}(y)\right).
\end{align}

We are ready to state the following Voronoi summation formula due to Miller--Schmid \cite{MS06} and Goldfeld--Li \cite{GL06}.

\begin{lemma}\label{lm:Voronoi3}
Let $g$ be a compactly supported smooth function on $\bR^+.$ For $(a,q)=1,$ we have
\begin{align}\label{eq:Voronoi3}
\sum_{m\geqslant1} \lambda_1(1,m)\ue\Big(\frac{am}{q}\Big)g(m)=&q\sum_{\pm}\sum_{m_1|q}\sum_{m_2\geqslant1} \frac{\lambda_1(m_2,m_1)}{m_1m_2}S(\bar a,\pm m_2; q/m_1)G_\pm\Big(\frac{m_1^2m_2}{q^3}\Big),
\end{align}
where $\bar{a}$ denotes the multiplicative inverse of $a\bmod{q}$ and 
\begin{align}\label{eq:Klsum}
S(m,n;c)=\sideset{}{^*}\sum_{x\bmod c}\ue\Big(\frac{mx+n\overline{x}}{c}\Big)\end{align}
denotes the classical Kloosterman sum.
\end{lemma}

Although the Ramanujan--Selberg conjecture is not proved yet for $\pi_1$, we have the following alternative estimates; see \cite[Remark 12.1.8]{Go06} or \cite[Theorem 2]{Mo02}.
\begin{lemma}\label{lm:second-GL3}
We have
$$
\sum_{n\leq N}|\lambda_1(n,1)|^2\ll N^{1+\varepsilon}
$$
and
$$
\sum_{n\leq N}|\lambda_1(1,n)|^2\ll N^{1+\varepsilon},
$$
where the implied constants depend on the form $\pi_1$ and $\varepsilon$.
\end{lemma}

We now turn to $SL(2,\bZ)$. For the sake of exposition we only present the case of Maass forms, and the case of holomorphic forms is just similar or even simpler. Furthermore, for technical simplicity, we only restrict to the case of full level. Let $\pi_2$ be a Maass cusp form with Laplace eigenvalue $\frac{1}{4}+t^2\geq 0$, and with Fourier expansion
$$
\sqrt{y}\sum_{n\neq 0}\lambda_2(n)K_{it}(2\pi|n|y)\ue(nx).
$$     
We will use the following Voronoi type summation formula (see Meurman \cite{Me88}).  
\begin{lemma}\label{lm:Voronoi2}
Let $h$ be compactly supported smooth function on $\bR^+$. For $(a,q)=1,$ we have
\begin{align}
\label{eq:Voronoi2}
\sum_{n=1}^\infty \lambda_2(n)\ue\Big(\frac{an}{q}\Big)h(n)=\frac{1}{q}\sum_{\pm}\sum_{n=1}^\infty \lambda_2(\mp n)\ue\Big(\pm\frac{\overline{a}n}{q}\Big)H^{\pm}\Big(\frac{n}{q^2}\Big)
\end{align}
where $\bar{a}$ is the multiplicative inverse of $a\bmod{q}$, and
\begin{align}
H^-(y)=&\frac{-\pi}{\cosh(\pi t)}\int_\bR h(x)\{Y_{2it}+Y_{-2it}\}\left(4\pi\sqrt{xy}\right)\ud x\label{eq:H-}\\
H^+(y)=&4\cosh(\pi t)\int_\bR h(x)K_{2it}\left(4\pi\sqrt{xy}\right)\ud x.\label{eq:H+}
\end{align}
\end{lemma}

\begin{remark}
If $g$ is supported in $[X,2X]$, satisfying $x^jg^{(j)}(x)\ll_j 1$, then the integral transform $G_\pm(y)$ satisfies
\begin{align}
\label{eq:Gtransform-bound}
G_{\pm}(y)\ll \sqrt{yX}(1+yX)^{-A}
\end{align}
for any fixed $A>0$ (see \cite[Remark 1]{Mu13a} for comments). Therefore,
the sums on the right hand side of \eqref{eq:Voronoi3} are essentially supported on $m_1^2m_2\ll q^3(qX)^{\varepsilon}/X$ (where the implied constant depends on $\pi_1$ and $\varepsilon$), and the contribution from the terms with $m_1^2m_2\gg q^3(qX)^{\varepsilon}/X$ is negligibly small, say $O((qX)^{-A})$ for any $A>0$.

If $h$ is supported in $[Y,2Y]$, satisfying $y^jh^{(j)}(y)\ll_j 1$, one has
\begin{align}\label{eq:Htransform-bound}
H^{\pm}(y)\ll Y(1+|y|Y)^{-A}
\end{align}
for any fixed $A>0$. Therefore, the sums on the right hand side of \eqref{eq:Voronoi2} are essentially supported on $n\ll q^2(qY)^{\varepsilon}/Y$ (where the implied constant depends on $\pi_2$ and $\varepsilon$) and the contribution from the terms with $n\gg q^2(qY)^{\varepsilon}/Y$ is also negligibly small.
\end{remark}

The following lemma characterizes the non-correlations of Fourier coefficients with additive characters.
\begin{lemma}\label{lm:Wilton-Miller}
Uniformly in $\alpha\in[0,1],$ we have
\begin{align}
\label{eq:Miller}
\sum_{n\leqslant X}\lambda_1(1,n)\ue(\alpha n)\ll X^{\frac{3}{4}+\varepsilon}
\end{align}
and
\begin{align}
\label{eq:Wilton}
\sum_{n\leqslant X}\lambda_2(n)\ue(\alpha n)\ll X^{\frac{1}{2}+\varepsilon},
\end{align}
where the implied constants depends on $\varepsilon$ and the form.
\end{lemma}

The inequality \eqref{eq:Wilton} is classical in the case of holomorphic forms of full level due to Wilton \cite{Wi29}; and the case of Maass forms can be proved in a similar way, which can be found, for instance, in \cite[Proposition 4]{HM06}.
The inequality \eqref{eq:Miller} is due to Miller \cite{Mil06} in a slightly more general setting.

\smallskip

\section{Trace functions and exponential sums}
\label{sec:trace-expsums}

This section is devoted to introduce the terminology and concepts on trace functions over finite fields, as well as the estimates on certain averages of them. More precisely, we shall state the arithmetic exponent pairs for {\it composite} trace functions developed in \cite{WX16}, which will be employed in the latter estimates for certain exponential sums.

\subsection{Trace functions}
Let $p$ be a prime and $\ell\neq  p$ an auxiliary prime, and fix an isomorphism $\iota : \overline{\bQ}_\ell\rightarrow\bC$.  Let $\cF$ be an $\ell$-adic middle-extension sheaf on $\bA^1_{\bF_p}$, which is pure of weight zero and of rank $\rank(\cF).$ The trace function associated to $\cF$ is defined to be
\begin{align*} K(x) := \iota((\tr\,\cF)(\bF_p, x))\end{align*}
for $x\in\bF_p$, following the manner of Katz \cite[Section 7.3.7]{Ka88}.
The (analytic) conductor of $\cF$, which is also called the conductor of $K$, is defined to be
\begin{align*}  
\fc(\cF) := \rank(\cF) + \sum_{x\in S(\cF)} (1+\Swan_x(\cF)),
\end{align*}
where $S(\cF)\subset\bP^1(\overline{\bF}_p)$ denotes the $($finite$)$ set of singularities of $\cF$, 
and $\Swan_x(\cF)$ $(\geqslant 0)$ denotes the Swan conductor of $\cF$ at $x$ $($see {\rm \cite{Ka80}}$).$

Let $q\geqslant3$ be a squarefree number. For each prime factor $p$ of $q$, one we may introduce an $\ell$-adic middle-extension sheaf $\cF_p$ on $\bA^1_{\bF_p}$, as well as its trace function $K_p(\cdot)$. The {\it composite} trace function $K=K_q(\cdot)$ is then defined by the product
\begin{align}\label{eq:def-Kq}
K_q(n)=\prod_{p\mid q}K_p(n),\end{align}
and conductor of $K$ can be defined in terms of the conductor of each $\cF_p$.

There are many typical examples of trace functions and regarding the applications in this paper, we only focus the one example: Kloosterman sums given by \eqref{eq:Klsum}. However, for the sake of geometric discussions, we will focus on the normalized sum
\begin{align*}\kl(n,q)=\frac{S(n,1;q)}{\sqrt{q}}\end{align*}
for squarefree $q\geqslant3$, which is a composite trace function mod $q$.
Recall the classical Weil bound
\begin{align*}|\kl(n,q)|\leqslant\tau(q).\end{align*}
More deeply, according to Deligne, there exists an $\ell$-adic middle-extension sheaf $\cK\ell$ modulo $p$, called a Kloosterman
sheaf with the corresponding trace function
\begin{align*}K_{\cK\ell}(x)=\kl(x,p),\ \ x\in\bF_p^\times.\end{align*}
Such a sheaf was constructed by Deligne \cite{De80}. According to him, $\cK\ell$ is geometrically irreducible, of rank $2$, and with conductor bounded by $5$.

\subsection{$\ell$-adic transforms and conductors}
For a non-trivial additive character $\psi$  and a function $f: \bF_p\rightarrow\bC$,
we define the Fourier transform $\ft_\psi(f):\bF_p\rightarrow\bC$ by
\begin{align*}\ft_\psi(f)(t)=\frac{-1}{\sqrt{p}}\sum_{x\in\bF_p}f(x)\psi(tx)\end{align*}
for $t\in\bF_p$.
According to Deligne \cite[3.4.1]{De80},  a middle-extension sheaf modulo $p$ of weight 0 is geometrically
 a direct sum of irreducible sheaves over $\bF_p$. Hence one can define a Fourier sheaf
modulo $p$ to be one where no such geometrically irreducible component is isomorphic to an Artin--Schreier sheaf $\cL_{\psi(aX)}$ for some $a\in\overline{\bF}_p.$
Such local Fourier transforms were studied in depth by Laumon \cite{La80},
Brylinski and Katz \cite{Ka80,Ka88}, and shown to satisfy the following properties (many of which
are, intuitively, analogues of classical properties of the Fourier transform):

\begin{lemma}\label{lm:fouriertransform}
Let $\psi$ be a non-trivial additive character of $\bF_p$ and $\cF$ a Fourier sheaf on $\bA_{\bF_p}^1.$  Then there exists an $\ell$-adic sheaf
$$\cG=\ft_\psi(\cF)$$
called the Fourier transform of $\cF$, which is also an $\ell$-adic Fourier sheaf, with
the property that 
$$K_{\ft_\psi(\cF)}(y)=\ft_\psi(K_\cF)(y)=\frac{1}{\sqrt{p}}\sum_{x\in\bF_p}K_\cF(x)\psi(yx).$$
Furthermore, we have

$(1)$ The sheaf $\cG$ is geometrically irreducible if and only if
$\cF$ is;

$(2)$ The Fourier transform is involutive, in the sense that we have a canonical arithmetic isomorphism
$$\ft_\psi(\cG)\simeq[\times(-1)]^*\cF,$$
where $[\times(-1)]^*$ denotes the pull-back by the map $x\mapsto-x;$

$(3)$ We have
\begin{align}\label{eq:cond-fouriertranform}
\fc(\ft_\psi(\cF))\leqslant10\fc(\cF)^2.\end{align}
\end{lemma}

\proof The last claim was proved by Fouvry, Kowalski and Michel \cite{FKM15} using the theory of local Fourier transforms developed by Laumon \cite{La80}, and the others have been established for instance in \cite[Theorem 8.4.1]{Ka88}.
\endproof

The inequality (\ref{eq:cond-fouriertranform}) is essential in analytic applications, since it implies that if $p$ varies but
$\cF$ has a bounded conductor, so do the Fourier transforms.

\subsection{Arithmetic exponent pairs for trace functions}
Let $q$ be a squarefree number and $K$ defined by \eqref{eq:def-Kq}. For a given interval $I$, we consider the following average of $K$ over $I$:
\begin{align*}
\fS(K_q,W_\delta)=\sum_{n\in I}K_q(n)W_\delta(n),
\end{align*}
where $W_\delta$ is an arbitrary function defined over $\bZ/\delta\bZ$ satisfying $\|W_\delta\|_\infty\leqslant1$. 
Roughly speaking, the P\'olya--Vinogradov bound is non-trivial for $|I|>q^{\frac{1}{2}+\varepsilon}$ (at least for $\delta=1$). The situation is much easier if $q$ allows certain factorizations. In \cite{WX16}, we developed the method of arithmetic exponent pairs for averages of composite trace functions that may go far beyond the P\'olya--Vinogradov bound, as long as the factorization of $q$ is good enough. Such observation can be at least dated back to Heath-Brown and his proof of Weyl-type subconvex bound for Dirichlet $L$-functions to well-factorable moduli.

Recall that a middle-extension sheaf $\cF_p$ on $\bA_{\bF_p}^1$,
which is pointwise pure of weight $0$ is said to be $d$-amiable if no geometrically irreducible component of $\cF_p$ is geometrically isomorphic to an Artin--Schreier sheaf of the form $\cL_{\psi(P)},$ where $P\in\bF_p[X]$ is a polynomial of degree $\leqslant d.$ In such case, we also say the associated trace function $K_p$ is $d$-amiable. The composite trace function $K_q$ is called to be $d$-amiable if $K_p$ is $d$-amiable for each $p\mid q.$ 
In addition, a sheaf $($or its associated trace function$)$ is said to be $\infty$-amiable if it is amiable for any fixed $d\geqslant1.$

We are ready to state the arithmetic exponent pairs developed in \cite{WX16}.
\begin{lemma}\label{lm:AEP}
Let $\eta>0$ be a sufficiently small number. Suppose $q$ is a squarefree number with no prime factors exceeding $q^\eta$, and $K=K_q(\cdot)$ is an $\infty$-amiable trace function $\bmod q$. For $|I|<q\delta,$ there exists $(\kappa,\lambda,\nu,\mu)$ such that
\begin{align}\label{eq:exponentpair}
\fS(K_q,W_\delta)\ll q^\varepsilon\Big(\frac{q}{|I|}\Big)^\kappa|I|^\lambda\delta^\nu\|\widehat{W}_\delta\|_\infty^\mu,
\end{align}
where $\varepsilon>0$ and the implied constant depends only on $\varepsilon$ and the conductor of $\cF_p$ for each $p\mid q$.

In particular, we may take $(\kappa,\lambda,\nu,\mu)=(\frac12,\frac12,\frac12,1),
(\frac{11}{30}, \frac{16}{30}, \frac{1}{6}, 1)$ and $(\frac{2}{18}, \frac{13}{18}, \frac{11}{28}, 0)$.
\end{lemma}

\begin{remark}
We call $(\kappa,\lambda,\nu,\mu)$ satisfying \eqref{eq:exponentpair} to be an arithmetic exponent pair for $(K_q,W_\delta)$. 
The values of $\nu,\mu$ are usually not too large so that these do not impact the applications. Note that one can take $(\kappa,\lambda,\nu,\mu)=(0,1,0,0)$ trivially and a sequence of exponent pairs can be produced starting from $(0,1,0,0)$ by virtue of $A$- and $B$-processes in the $q$-analogue of the van der Corput method. In the following application, however, we only invoke the choice $(\kappa,\lambda,\nu,\mu)=(\frac12,\frac12,\frac12,1).$

Since we do not input extra conditions on $W_\delta$, it is not a bad choice to utilize the trivial bound $\|\widehat{W}_\delta\|_\infty\leqslant\sqrt{\delta}.$
\end{remark}

\subsection{Babies of Kloosterman sums}

By virtue of Kloosterman sums, we define two relevant algebraic sums:
\begin{itemize}
\item For $d\mid c,$ define
\begin{align}\label{eq:S(h,n,m;c,d)}
\sS(h,n,m;c,d)=\sideset{}{^*}\sum_{x\bmod c}S(m,x;d)\ue\Big(\frac{h\overline{x}-nx}{c}\Big). 
\end{align}

\item \begin{align*}
T(a,b,m;c)&=\frac{1}{c}~\sideset{}{^*}\sum_{x\bmod c}S(\overline{x}+a,-b;c)\ue\Big(\frac{-mx}{c}\Big).
\end{align*}
\end{itemize}

For $d\mid c,$ put $\ell=c/d$. If $(d,\ell)=1,$ from the Chinese remainder theorem it follows that
\begin{align*}
\sS(h,n,m;c,d)=S(h,-n\overline{d}^2;\ell)\sideset{}{^*}\sum_{x\bmod d}S(m,x;d)\ue\Big(\frac{\overline{\ell}(h\overline{x}-nx)}{d}\Big). 
\end{align*}
Furthermore, opening Kloosterman sums and rearranging the summations, we may conclude the following lemma.
\begin{lemma}
Let $c=d\ell$ with $(d,\ell)=1$, we have
\begin{align*}
\sS(h,n,m;c,d)=d\cdot S(h,-n\overline{d}^2;\ell)T(n\overline{\ell},h\overline{\ell},m;d).\end{align*}
\end{lemma}

The Chinese remainder theorem also yields the following twisted multiplicativity of 
$T(a,b,m;c).$
\begin{lemma}
For $c=c_1c_2$ with $(c_1,c_2)=1,$ we have
\begin{align*}
T(a,b,m;c_1c_2)&=T(a,b\overline{c}_2^2,m\overline{c}_2;c_1)T(a,b\overline{c}_1^2,m\overline{c}_1;c_2).
\end{align*}
\end{lemma}

For the sake of subsequent applications of arithmetic exponent pairs, we would like to determine when $T(a_1,b_1,y;p)$ and $T(a_2,b_2,y;p)$, as functions in $x$, do not correlate for given tuples $(a_1,b_1)$ and $(a_2,b_2).$

Note that $y\mapsto T(a_1,b_1,y;p)$ is Fourier transform of $S(\overline{x}+a,-b;p)/\sqrt{p}$. For $p\nmid b$, we have
\[\frac{S(\overline{x}+a,-b;p)}{\sqrt{p}}=\kl(-b(\overline{x}+a),p),\]
which is a trace function of $[\gamma\cdot x]^*\cK\ell$ with 
\begin{align*}
\gamma=
\begin{pmatrix}
ab & b\\
-1 & 0
\end{pmatrix}.
\end{align*}
In view of Lemma \ref{lm:fouriertransform}, it suffices to consider $[\gamma\cdot x]^*\cK\ell$.

On the other hand, it was proved by Katz that
there does not exist a rank 1 sheaf $\cL$ and a geometric isomorphism
\begin{align*}
[\gamma\cdot x]^*\cK\ell\simeq\cK\ell\otimes\cL
\end{align*}
for $\mathbf{1}\neq\gamma\in PGL_2(\bF_p).$ As a consequence, if $\gamma\neq \mathbf{1}$, we find,  for any rank 1 sheaf $\cL,$ the triple tensor sheaf
\begin{align*}
[\gamma\cdot x]^*\cK\ell\otimes\cK\ell\otimes\cL
\end{align*}
is geometrically irreducible, of rank $\geqslant2,$ and thus $\infty$-amiable.

Put
\begin{align*}
\gamma_1=
\begin{pmatrix}
a_1b_1 & b_1\\
-1 & 0
\end{pmatrix},\ \ \ 
\gamma_2=
\begin{pmatrix}
a_2b_2 & b_2\\
-1 & 0
\end{pmatrix}
\end{align*}
and suppose $p\nmid b_1b_2.$ Hence $\gamma_1,\gamma_2\in PGL_2(\bF_p)$.
We would like to determine when
\begin{align*}
[\gamma_1\cdot x]^*\cK\ell\otimes[\gamma_2\cdot x]^*\cK\ell\otimes\cL
\end{align*}
is $\infty$-amiable. After a bijective change of variable, it suffices to check $\gamma_2\gamma_1^{-1}$ is the identity or not.
In fact,
\begin{align*}
\gamma_2\gamma_1^{-1}=
\begin{pmatrix}
a_2b_2 & b_2\\
-1 & 0
\end{pmatrix}
\begin{pmatrix}
0 & -1\\
1/b_1 & a_1
\end{pmatrix}
=\begin{pmatrix}
b_2/b_1 & b_2(a_1-a_2)\\
0 & 1
\end{pmatrix},
\end{align*}
which is identity if and only if $a_1\equiv a_2\bmod p$ and $b_1\equiv b_2\bmod p$.
Hence we may conclude the following lemma.
\begin{lemma}\label{lm:amiable}
Suppose $p\nmid b_1b_2.$
If $a_1\not\equiv a_2\bmod p$ or $b_1\not\equiv b_2\bmod p$,
the trace function
\begin{align*}
y\mapsto T(a_1,b_1,y;p)\overline{T(a_2,b_2,y;p)}
\end{align*}
is $\infty$-amiable.
\end{lemma}

\begin{remark}
Suppose $p\nmid h\ell k_1k_2$. Since $T(n_j\overline{\ell},h\overline{\ell k_j^2},y\overline{k}_j;p)=T(n_jk_j\overline{\ell},h\overline{\ell k_j^3},y;p)$ 
for $j=1,2$. It follows from Lemma \ref{lm:amiable} that
$y\mapsto T(n_1\overline{\ell},h\overline{\ell k_1^2},y\overline{k}_1;p)
\overline{T(n_2\overline{\ell},h\overline{\ell k_2^2},y\overline{k}_2;p)}$
is $\infty$-amiable if
\[n_1\equiv n_2\bmod p,\ \ k_1^3\equiv k_2^3\bmod p.\]
In particular, if $p\equiv2\bmod3$, the later congruence condition is equivalent to
$k_1\equiv k_2\bmod p.$ This observation will be invoked in the proof Theorem \ref{thm:individual}. 
\end{remark}

\smallskip
\section{Applying circle method}
\label{sec:circlemethod}

\subsection{Jutila's variant of circle method}
In the study of $\sD_h(X)$, we would like to detect the condition $m=n\in\bZ$ starting from the trivial identity
\begin{align*}
\int_0^1\ue(\alpha(m-n))\ud\alpha=
\begin{cases}
1,\ \ &\text{if }m=n,\\
0,&\text{if }m\neq n.
\end{cases}
\end{align*}
The circle method is devoted to the decomposition of the unit interval $[0,1]$ in a certain way such that the subsequent evaluations can be proceeded non-trivially. There is a flexible invariant due to Jutila \cite{Ju92,Ju97} with overlapping intervals, which has found many important applications in the analytic theory of automorphic forms.

Let $Q\geqslant1$. For any moduli set $\sQ \subseteq [Q,2Q]$ and a positive real number $\varDelta$ with $ Q^{-2}\ll\varDelta\ll Q^{-1}$, we define the function
$$
I_{\sQ,\varDelta} (x)=\frac{1}{2\varDelta\varPhi}\sum_{q\in\sQ}\;\sideset{}{^*}\sum_{a\bmod{q}}\mathbf{1}_\varDelta\Big(\frac{a}{q}-x\Big),
$$
where $\varPhi=\sum_{q\in\sQ}\varphi(q)$ and $\mathbf{1}_\varDelta$ is the characteristic function of $[-\varDelta,\varDelta]$.

If we view the fractions $a/q$ as random variables, the expected deviation of $I_{\sQ,\varDelta} (x)$ from 1 is $O(1/(\varDelta Q^2))$. In fact, Jutila proved this is the case in the following average sense, if $\varPhi\gg Q^{2-\varepsilon}$, which holds upon many interesting choices.
The proof can be found in \cite{Ju92}.

\begin{lemma}\label{lm:jutila}
Let $Q^{-2}\ll \varDelta\ll Q^{-1} $. Then we have
\begin{align}
\label{variance}
\frac{1}{\varDelta Q^2}\int_0^1\left|1-I_{\sQ,\varDelta}(x)\right|^2\ud x\llcurly \frac{1}{(\varDelta\varPhi)^2}.
\end{align}
\end{lemma}

\subsection{Setting up the circle method}

Let $W$ be a smooth function supported in $[\frac{1}{2},3]$ satisfying $W(x)=1$ for $x\in [\frac{2}{3},\frac{5}{2}]$. Then we may write
\begin{align}\label{eq:unit-identity}
\sD_h(X)
&=\int_0^1\ue(xh)\fS_1(x,X)\fS_2(x,X)\ud x,
\end{align}
where
\begin{align}\label{eq:S1(x,X)}
\fS_1(x,X)=\sum_{m\geqslant1}\lambda_1(1,m)\ue(xm)V\left(\frac{m}{X}\right)
\end{align}
and
\begin{align}\label{eq:S2(x,X)}
\fS_2(x,X)=\sum_{n\geqslant1}\lambda_2(n)\ue(-xn)W\left(\frac{n}{X}\right).
\end{align}

In \cite{Mu13a}, Munshi constructed $\sQ$ to be a set of products of distinct primes with prescribed sizes. In order to figure a general framework, we here allow $\sQ$ to be a collection of squarefree numbers contained in $[Q/2,Q]$. Let
$\varDelta>0$ be parameter to be chosen later such that $Q^{-2}\ll \varDelta\ll Q^{-1}.$
Define
\begin{align*}
\sD_h^*(X)
&=\int_0^1I_{\sQ,\varDelta}(x)\ue(xh)\fS_1(x,X)\fS_2(x,X)\ud x.
\end{align*}

We now first consider the difference between $\sD_h(X)$ and $\sD_h^*(X).$
From \eqref{eq:Wilton} and partial summation, it follows that
$$\fS_2(x,X)\llcurly X^{\frac{1}{2}},
$$
which is uniform in $x\in[0,1]$. We then have
\begin{align*}
|\sD_h(X)-\sD^*_h(X)|\llcurly X^{\frac{1}{2}}\int_0^1|\fS_1(x,X)|\left|1- I_{\sQ,\varDelta}(x)\right|\ud x.
\end{align*}
By Cauchy's inequality, it follows from Lemmas \ref{lm:jutila} and \ref{lm:second-GL3} that
\begin{align*}
|\sD_h(X)-\sD^*_h(X)|^2
&\llcurly \frac{(QX)^2}{\varDelta\varPhi^2}.
\end{align*}

If we specialize the moduli set $\sQ$ such that $\varPhi\gg Q^{2-\varepsilon}$, we then conclude the following approximation.
\begin{lemma}
\label{lm:Dh(X)-approximation}
For $Q^{-2}\ll \varDelta\ll Q^{-1}$ and $\varPhi\gg Q^{2-\varepsilon}$, we have
\begin{align}\label{eq:Dh(X)-approximation}
\sD_h(X)=\sD_h^*(X)+O\Big(\frac{X^{1+\varepsilon}}{\sqrt{\varDelta}Q}\Big),
\end{align}
where the $O$-constant depends on $\varepsilon.$
\end{lemma}

\begin{remark}
We will specialize $\varDelta=1/X$, which implies an error term $O(Q^{-1}X^{\frac{3}{2}+\varepsilon})$ in \eqref{eq:Dh(X)-approximation}. To beat the trivial bound of $\sD_h(X)$, we require that $Q\gg \sqrt{X}$, which is always kept in mind henceforth.
\end{remark}

\smallskip
\section{Transformation of $\sD_h^*(X)$}
\label{sec:transformation}

\subsection{Applying Voronoi summation formulae}
In this section, we make some initial transformations of
$\sD_h^*(X)$ by appealing to the Voronoi summation formulae on $GL(3)$
and $GL(2)$. 

From the definition of $I_{\sQ,\varDelta}(x)$, we get
\begin{align*}
\sD_h^*(X)&=\frac{1}{2\varDelta L}\int_{-\varDelta}^{\varDelta}\sum_{q\in\sQ}~\sideset{}{^*}\sum_{a\bmod q}\ue\Big(\frac{ah}{q}\Big)\mathop{\sum\sum}_{m,n\in\bZ}\lambda_1(1,m)\lambda_2(n)\ue\Big(\frac{a(m-n)}{q}\Big)\\
&\ \ \ \ \ \times V\left(\frac{m}{X}\right)W\left(\frac{n}{Y}\right)\ue(\alpha(m-n+h))\ud\alpha.\end{align*}
For any fixed $\alpha\in[-\varDelta,\varDelta],$ put
\begin{align}\label{eq:D(h,alpha)}
\sD_{h,\alpha}^*(X)&=\sum_{q\in\sQ}~\sideset{}{^*}\sum_{a\bmod q}\ue\Big(\frac{ah}{q}\Big)\mathop{\sum\sum}_{m,n\in\bZ}\lambda_1(1,m)\lambda_2(n)\ue\Big(\frac{a(m-n)}{q}\Big)g(m)h(n)\end{align}
with
\begin{align}\label{eq:choices:gh}
g(x)=V\left(\frac{x}{X}\right)\ue(\alpha x),\ \ h(y)=W\left(\frac{y}{Y}\right)\ue(-\alpha y),\end{align}
so that
\begin{align}
\label{eq:Dh*}
\sD_h^*(X)&=\frac{1}{2\varDelta L}\int_{-\varDelta}^{\varDelta}\ue(\alpha h)\sD_{h,\alpha}^*(X)\ud\alpha.\end{align}
We are now in a position to apply Voronoi summation formulae (Lemmas \ref{lm:Voronoi3} and \ref{lm:Voronoi2}) to sums over $m,n$ in \eqref{eq:D(h,alpha)}, getting
\begin{align*}
\sD_{h,\alpha}^*(X)&=\mathop{\sum\sum}_{\sigma_1,\sigma_2\in\{+,-\}}\sum_{q\in\sQ}\sum_{m_1|q}\sum_{m_2\geqslant1} \frac{\lambda_1(m_2,m_1)}{m_1m_2}\\
&\ \ \ \ \ \times\sum_{n\geqslant1} \lambda_2(-\sigma_2n)\sS(h,n,m_2;q,q/m_1)G_{\sigma_1}\Big(\frac{m_1^2m_2}{q^3}\Big)H^{\sigma_2}\Big(\frac{n}{q^2}\Big),
\end{align*}
where $G_\pm$ and $H^\pm$ are integral transforms defined by \eqref{eq:Gpm}, \eqref{eq:H-} and \eqref{eq:H+} with respect to the choices in \eqref{eq:choices:gh}, and the exponential sum $\sS(h,n,m_2;q,q/m_1)$ is given by \eqref{eq:S(h,n,m;c,d)}.

\subsection{Well-factorable moduli coming into play}
Due to the choices of $\sigma_1,\sigma_2$ in the above expression for $\sD_{h,\alpha}^*(X)$, there will be four sub-contributions with strong similarities; we choose only one to ease the complexity of presentation. More precisely, we consider 
\begin{align*}
\varSigma_{h,\alpha}(X)&=\sum_{q\in\sQ}\sum_{\ell|q}\mathop{\sum\sum}_{m,n\geqslant1}\frac{\lambda_1(m,\ell)\lambda_2(n)}{m\ell}\sS(h,n,m;q,q/\ell)G_+\Big(\frac{\ell^2m}{q^3}\Big)H^-\Big(\frac{n}{q^2}\Big).\end{align*}
Alternatively, we have the expression
\begin{align*}
\varSigma_{h,\alpha}(X)&=\sum_{\ell d\in\sQ}\mathop{\sum\sum}_{m,n\geqslant1}\frac{\lambda_1(m,\ell)\lambda_2(n)}{m\ell}\sS(h,n,m;\ell d,d)G_+\Big(\frac{m}{\ell d^3}\Big)H^-\Big(\frac{n}{\ell^2d^2}\Big).\end{align*}

The moduli set $\sQ$ is chosen to be the set of squarefree numbers in $(Q/2,Q],$
whose prime factors are $\equiv2\bmod3$ and not exceeding $Q^\eta$
for any small $\eta>0.$ In particular, we may take $\eta=\varepsilon^2.$ It is known that
\[\varPhi:=\sum_{q\in\sQ}\varphi(q)\gg Q^2\]
with an implied constant depends only on $\eta$.
 By dyadic devices, we consider
\begin{equation}\label{eq:Sigma_h,alpha(X;M,N)}
\begin{split}
\varSigma_{h,\alpha}(X;M,N)&=\sum_{\ell\sim L}\sum_{d\sim D}\sum_{k\sim K}\sum_{m\sim M}\sum_{n\sim N}\frac{\lambda_1(m,\ell)\lambda_2(n)}{m\ell}\sS(h,n,m;\ell dk,dk)\\
&\ \ \ \ \ \ \times G_+\Big(\frac{m}{\ell d^3k^3}\Big)H^-\Big(\frac{n}{\ell^2d^2k^2}\Big)\end{split}
\end{equation}
for all $L\leqslant Q$ with $LDK\asymp Q$, where $D,K$ can be specialized on our demand. In addition, the sums are void unless $\mu^2(\ell dk)=1,$
which we always assume henceforth.
From the supports and decays of $G_+,H^-$ given by \eqref{eq:Gtransform-bound} and \eqref{eq:Htransform-bound}, we may assume
\begin{align}\label{eq:MN-size}
M\llcurly X^{-1}QD^2K^2,\ \ N\llcurly X^{-1}Q^2.
\end{align}

Our next goal is to bound $\varSigma_{h,\alpha}(X;M,N)$
from above uniformly in $\alpha\in[-\varDelta,\varDelta]$. Collecting all possible tuples $(M,N)$ and integrating over $\alpha$ trivially, one can obtain an upper bound for $\sD^*_h(X)$, from which and Lemma \ref{lm:Dh(X)-approximation} the theorem follows by optimizing the parameters.

\smallskip

\section{Proof of Theorem \ref{thm:individual}}
\subsection{Invoking bilinear structures}
By Cauchy's inequality, it follows from \eqref{eq:Sigma_h,alpha(X;M,N)} that
\begin{align*}
|\varSigma_{h,\alpha}(X;M,N)|^2\leqslant\varSigma_1\varSigma_2,\end{align*}
where
\begin{align*}
\varSigma_1&=\sum_{\ell\sim L}\sum_{d\sim D}\sum_{m\sim M}\left|\frac{\lambda_1(m,\ell)}{m\ell}\right|^2\end{align*}
and
\begin{align*}
\varSigma_2&=\sum_{\ell\sim L}\sum_{d\sim D}\sum_{m\sim M}\Bigg|\sum_{k\sim K}\sum_{n\sim N}\lambda_2(n)\sS(h,n,m;\ell dk,dk)G_+\Big(\frac{m}{\ell d^3k^3}\Big)H^-\Big(\frac{n}{\ell^2d^2k^2}\Big)\Bigg|^2.\end{align*}

On one hand, it follows from Lemma \ref{lm:second-GL3} that
\begin{align}\label{eq:Sigma1-upperbound}
\varSigma_1
&\llcurly\frac{D}{(LM)^2}\sum_{d}\sum_{\ell\sim L/d}\sum_{m\sim M/d}|\lambda_1(m,1)|^2|\lambda_1(\ell,1)|^2\llcurly\frac{D}{LM}.
\end{align}
Squaring out and switching summations, 
\begin{align}\label{eq:Sigma2-upperbound}
\varSigma_2&\llcurly X^2\sum_{\ell\sim L}\sum_{d\sim D}\mathop{\sum\sum}_{k_1,k_2\sim K}\mathop{\sum\sum}_{n_1,n_2\sim N}|\lambda_2(n_1)||\lambda_2(n_2)||\sB(\ell,d,\bk,\bn)|,\end{align}
where
\begin{align*}
\sB(\ell,d,\bk,\bn)&=\sum_{m\sim M}\sS(h,n_1,m;\ell dk_1,dk_1)\overline{\sS(h,n_2,m;\ell dk_2,dk_2)}G_+\Big(\frac{m}{\ell d^3k_1^3}\Big)\overline{G_+\Big(\frac{m}{\ell d^3k_2^3}\Big)}.\end{align*}

\subsection{Estimates for $\varSigma_2$ and $\sD^*_h(X)$}
The next step is to capture the oscillations within $\sB(\ell,d,\bk,\bn)$ by virtue of the method of arithmetic exponent pairs given in Lemma \ref{lm:AEP}.
Put $k_0=(k_1,k_2)$ and $k_1=k_0k_1',k_2=k_0k_2',$ so that $(k_1',k_2')=1.$
We then decompose 
\begin{align*}
\sS(h,n_1,m;\ell dk_1,dk_1)
&=dk_1\cdot S(h,-n_1\overline{(dk_1)^2};\ell)T(n_1\overline{\ell},h\overline{\ell},m;dk_1)\\
&=dk_1\cdot S(h,-n_1\overline{(dk_1)^2};\ell)T(n_1\overline{\ell},h\overline{\ell k_1^2},m\overline{k}_1;d)\\
&\ \ \ \ \ \cdot
T(n_1\overline{\ell},h\overline{\ell(dk_1')^2},m\overline{dk_1'};k_0)
T(n_1\overline{\ell},h\overline{\ell(dk_0)^2},m\overline{dk}_0;k_1'),\end{align*}
and we have a relevant decomposition for $\sS(h,n_2,m;\ell dk_2,dk_2)$.
As a function in $m$, the product $\sS(h,n_1,m;\ell dk_1,dk_1)\overline{\sS(h,n_2,m;\ell dk_2,dk_2)}$ is well-defined over $\bZ/dk\bZ$ with $k=k_0k_1'k_2'.$

Due to the possible correlation between $T(n_1\overline{\ell},h\overline{\ell k_1^2},m\overline{k}_1;d)$ and $T(n_2\overline{\ell},h\overline{\ell k_2^2},m\overline{k}_2;d),$ 
we would like to divide the tuples $(n_1,n_2,k_1,k_2)$ into the following three cases:
\begin{itemize}
\item $n_1=n_2,k_1=k_2,$
\item $n_1\neq n_2,k_1=k_2,$
\item $k_1\neq k_2.$
\end{itemize}
Correspondingly, we denote by $\varSigma_{21},\varSigma_{22}$ and $\varSigma_{23}$ the relevant contributions from these tuples to the RHS in \eqref{eq:Sigma2-upperbound}.

In Case I, we have $k_0=k$ and $k_1'=k_2'=1$, and we invoke the trivial bound
\begin{align*}
|\sS(\cdots)\overline{\sS(\cdots)}|&\llcurly\ell d^2k_1k_2\cdot(h,n_1,\ell),
\end{align*}
getting
\begin{align*}
\varSigma_{21}&\llcurly M^2NLX^3.\end{align*}

In Case II, we also have $k_0=k=k_1=k_2$ and $k_1'=k_2'=1$. Put $\fp=(dk,n_1-n_2)$ and $\fq=dk/\fp.$ In such case, we may have $(n_1-n_2,\fq)=1$ and write
\begin{align*}
&\sS(h,n_1,m;\ell dk_1,dk_1)\overline{\sS(h,n_2,m;\ell dk_2,dk_2)}\\
=&(dk)^2\cdot S(h,-n_1\overline{(dk)^2};\ell)S(h,-n_2\overline{(dk)^2};\ell)\\
&\ \ \ \ \times
T(n_1\overline{\ell},h\overline{\ell},m;dk)\overline{T(n_2\overline{\ell},h\overline{\ell},m;dk)}\\
=&(dk)^2\cdot S(h,-n_1\overline{(dk)^2};\ell)S(h,-n_2\overline{(dk)^2};\ell)|T(n_1\overline{\ell},h\overline{\ell\fq^2},m\overline{\fq};\fp)|^2\\
&\ \ \ \ \times T(n_1\overline{\ell},h\overline{\ell\fp^2},m\overline{\fp};\fq)\overline{T(n_2\overline{\ell},h\overline{\ell\fp^2},m\overline{\fp};\fq)}.\end{align*}
Thanks to Lemma \ref{lm:amiable}, we may apply Lemma \ref{lm:AEP}  with
\begin{align*}
q=\fq,\ \ \delta=\fp, \ \ (\kappa,\lambda,\nu,\mu)=(\frac12,\frac12,\frac12,1),
\end{align*}
\begin{align*}
K_q:x\mapsto T(n_1\overline{\ell},h\overline{\ell\fp^2},m\overline{\fp};\fq)\overline{T(n_2\overline{\ell},h\overline{\ell\fp^2},m\overline{\fp};\fq)}
\end{align*}
and
\begin{align*}
W_\delta:x\mapsto|T(n_1\overline{\ell},h\overline{\ell\fq^2},m\overline{\fq};\fp)|^2,
\end{align*}
getting
\begin{align*}
\sB_h(\ell,d,\bk,\bn)&\llcurly\frac{MX}{L(DK)^3}L(DK)^2(h,n_1,\ell)^{\frac{1}{2}}(h,n_2,\ell)^{\frac{1}{2}}\Big(\fq^{\frac{1}{2}}\fp+\frac{M}{\sqrt{DK}}\Big)\\
&\llcurly\frac{MX}{DK}(h,\ell)(dk)^{\frac{1}{2}}(dk,n_1-n_2)^{\frac{1}{2}}+\frac{M^2X}{(DK)^{\frac{3}{2}}}(h,\ell),\end{align*}
from which we conclude
\begin{align*}
\varSigma_{22}&\llcurly L(DK)^{\frac{1}{2}} MN^2X^3+\frac{L(MN)^2X^3}{\sqrt{DK}}.\end{align*}

In Case III, we can follow the arguments in Case II by pulling out $(d,k_1-k_2)$, the g.c.d. of $d$ and $k_1-k_2$, then apply Lemmas \ref{lm:AEP} and \ref{lm:amiable} with
\[q=k_1'k_2'd/(d,k_1-k_2),\ \ \delta=k_0(d,k_1-k_2),\ \ (\kappa,\lambda,\nu,\mu)=(\frac12,\frac12,\frac12,1),\]
and the relevant $K_q,W_\delta$, of which we do not intend to display the explicit shapes, getting
\begin{align*}
\sB_h(\ell,d,\bk,\bn)
&\llc\frac{MX}{DK}(h,\ell)\Big((k_1'k_2'd)^{\frac{1}{2}}(d,k_1-k_2)^{\frac{1}{2}}k_0+\frac{M}{\sqrt{dk_1'k_2'k_0}}\Big)\\
&\ll\frac{MX}{\sqrt{D}}(h,\ell)(d,k_1-k_2)^{\frac{1}{2}}+\frac{M^2X}{D^{\frac{3}{2}}K^2}(h,\ell)(k_1,k_2)^\frac{1}{2}.\end{align*}
Summing over $\ell,d,k_1,k_2,n_1,n_2$ in Case III, we derive that
\begin{align*}
\varSigma_{23}
&\llc LD^{\frac{1}{2}}K^2MN^2X^3
+\frac{L(MN)^2X^3}{\sqrt{D}}.\end{align*}

Collecting all above estimates for $\varSigma_{21},\varSigma_{22}$ and $\varSigma_{23}$, we find
\begin{align*}
\varSigma_2&\llc LM^2NX^3+LD^{\frac{1}{2}}K^2MN^2X^3
+\frac{L(MN)^2X^3}{\sqrt{D}},\end{align*}
from which and \eqref{eq:Sigma1-upperbound}, it follows that
\begin{align*}
|\varSigma_{h,\alpha}(X;M,N)|^2
&\llc X^3(DMN+D^{\frac{3}{2}} K^2N^2
+\sqrt{D}MN^2).\end{align*}
Collecting all possible tuples $(M,N)$ subject to \eqref{eq:MN-size}, we find
\begin{align*}
\varSigma_{h,\alpha}(X)&\llc(DQ^5X+D^{-\kappa}Q^{5+3\lambda-\kappa}X^{1+\kappa-\lambda}
+\sqrt{D}Q^7)^{\frac{1}{2}},\end{align*}
and the same upper bound also works for $\sD^*_{h,\alpha}(X).$
Integrating over $\alpha\in[-\varDelta,\varDelta],$
we obtain
\begin{align}\label{eq:D*h(X)-estimate}
\sD^*_h(X)&\llc(DQX+D^{-\frac{1}{2}}Q^2X
+\sqrt{D}Q^3)^{\frac{1}{2}}.\end{align}

\subsection{Concluding Theorem \ref{thm:individual}}

Taking $D=Q^{\frac{2}{3}}$ in \eqref{eq:D*h(X)-estimate} to balance the first two terms on the RHS, we get
\begin{align*}
\sD^*_h(X)& \llcurly Q^{\frac{5}{6}}X^{\frac{1}{2}}+Q^{\frac{5}{3}},
\end{align*}
from which and \eqref{eq:Dh(X)-approximation} we conclude
\begin{align*}
\sD_h(X)&\llcurly Q^{\frac{5}{6}}X^{\frac{1}{2}}+Q^{\frac{5}{3}}+\frac{X}{\sqrt{\varDelta}Q}.
\end{align*}

Taking $Q=X^{\frac{6}{11}}$ and $\varDelta=X^{-1}\geqslant Q^{-2}$, we arrive at
\begin{align*}
\sD_h(X)&\llcurly X^{\frac{21}{22}}
\end{align*}
as stated in Theorem \ref{thm:individual}.

\smallskip

\section{Remarks}
Regarding the estimates for exponential sums, Munshi \cite{Mu13a} appealed to the works of Adolphson--Sperber \cite{AS89} and Bombieri--Sperber \cite{BS95} on multi-dimensional exponential sums over finite fields. Alternatively, we utilize the arithmetic exponent pairs (Lemma \ref{lm:AEP}) developed in \cite{WX16}, which
allow us to deal with these exponential sums simultaneously. On the other hand, that the function in Lemma \ref{lm:amiable} is $\infty$-amiable is much more than what we need in the proof of Theorem \ref{thm:individual}, since we only utilize the exponent pair $(\frac12,\frac12,\frac12,1).$

As in the case of $GL(2)$, one can also consider the average of $\sD_h(X)$ over $h$. In fact, we may write
\begin{align*}
\sum_{h\in\bZ}|\sD_h(X)|^2&=\sum_{h\in\bZ}\left|\sum_{\substack{m,n\geqslant1\\m+h=n}}\lambda_1(1,m)\lambda_2(n)V\left(\frac{m}{X}\right)
W\left(\frac{n}{X}\right)\right|^2,
\end{align*}
where $W$ is a smooth function supported in $[\frac{1}{2},3]$ satisfying $W(x)=1$ for $x\in [\frac{2}{3},\frac{5}{2}]$.  From the orthogonality of additive characters, it follows that
\begin{align*}
\sum_{h\in\bZ}|\sD_h(X)|^2&=\int_0^1|\fS_1(x,X)\fS_2(x,X)|^2\ud x,
\end{align*}
where $\fS_1(x,X)$ and $\fS_2(x,X)$ are defined by \eqref{eq:S1(x,X)} and \eqref{eq:S2(x,X)}. We then conclude from Lemmas \ref{lm:Wilton-Miller} and \ref{lm:second-GL3} that
\begin{align*}
\sum_{h\in\bZ}|\sD_h(X)|^2&\llc X\int_0^1|\fS_1(x,X)|^2\ud x\llc X^2.
\end{align*}
This gives the square-root cancellation in $\sD_h(X)$ on average. As an immediate consequence, we also have
\begin{align}\label{eq:Dh(X)-average}
\sum_{h}\gamma(h)\sD_h(X)\llc X\Big(\sum_{h}|\gamma(h)|^2\Big)^{\frac12}\end{align}
for an arbitrary coefficient $\boldsymbol\gamma=(\gamma(h))\in\ell^2(\bZ)$.

It is expected that the approach in this paper can be utilized to improve \eqref{eq:Dh(X)-average} in some special cases. Unfortunately, 
a lot of information cannot be saved in the application of Lemma \ref{lm:Dh(X)-approximation}, although one can choose better exponent pairs while applying Lemma \ref{lm:AEP}. 

In fact, in the contexts of $GL(2)$, one can choose the moduli set to be consecutive integers (with certain necessary divisibility conditions at most), and one will encounter sums of Kloosterman sums after Voronoi, so that Kuznetsov trace formula is applicable since one is summing over consecutive integers. Hence, many more cancellations become possible thanks to Kloostermania. See \cite{Bl05,BHM07} for instances.

On the other hand, in the case of $GL(3)\times GL(2)$, Sun \cite{Su17} proved that
\begin{align}\label{eq:Dh(X)-smoothaverage}
\sum_{h\geqslant1}U\Big(\frac{h}{H}\Big)\sD_h(X)\ll X^{-A}
\end{align}
for any $A>0$, provided that $H>X^{\frac{1}{2}+\varepsilon}$, where $U$ is a fixed smooth function with compact support in $\bR^+.$ Instead of using Jutila's variant of the circle method, she employed Kloosterman's circle method in the version of Heath-Brown \cite{HB96}, in which case no approximation such as Lemma \ref{lm:Dh(X)-approximation} is required. The saving in \eqref{eq:Dh(X)-smoothaverage} comes from the estimate
\[\sum_{n\geqslant1}\lambda_2(n)V\Big(\frac{n}{X}\Big)\ll X^{-B},\ \ B>0\]
due to Booker \cite{Bo05}.

\end{document}